\documentclass[a4paper,12pt]{article}

\usepackage{amsmath}
\usepackage{amsfonts}
\usepackage{amssymb}

\begin{document}


\begin{center}

{ \Large \ttfamily
On Averaging for Hamiltonian Systems \\
with One Fast Phase and Small Amplitudes.
}
\vspace{5mm}

Jochen Br\"uning \\
Institut f\"ur Mathematik der
Humboldt--Universit\"at zu Berlin\\
Rudower Chausee 25, WBC(I.313) 12489 Berlin-Adlershof, Germany \\
E-mail: bruening{@}spectrum.mathematik.hu-berlin.de\\
\vspace{5mm}

Serguei Dobrokhotov \\
Institute for Problems in Mechanics of the
Russian Academy of Sciences,\\
pr.~Vernadskogo 101, 117526 Moscow,
Russia;\\
E-mail: dobr{@}ipmnet.ru
\vspace{5mm}

Michael Poteryakhin \\
Russian Research Center ''Kurchatov Institute'', \\
Kurchatov Square, 123182 Moscow, Russia;\\
E-mail: stpma{@}inse.kiae.ru

\end{center}

\section{ Introduction. }

   The problem of averaging for systems with one fast
   phase was considered from various points of view
   in many papers. The averaging method of Krylov and Bogolyubov
   \cite{BM} and methods of KAM theory originated this line of research,
   the most complete results were obtained by Neishtadt \cite{Ne}, where the
   coefficients are assumed real analytic. However, in many problems which
   are interesting from the point of view of applications,
   analytic dependence ceases to hold in the neighborhoods of
   some points.

   For instance, consider the motion of a particle under the influence of a
   small periodic electric and constant magnetic field \cite{BD}.
   The Hamiltonian of such a system has the form~:
   \begin{align*}
      H= \frac{1}{2}\left( (p_1+x_2)^2+ p_2^2 \right) +
         \varepsilon V(x_1,x_2).
   \end{align*}
   The canonical change of variables $ x_1=Q+y_1$, $ p_1=-y_2 $,
   $ x_2=P+y_2$, $ p_2=-Q $, $ Q= \sqrt{2I} \cos \varphi $,
   $ P=\sqrt{2I} \sin \varphi $, transforms it into:
   \begin{align*}
      H= I + \varepsilon V \left( \sqrt{2I} \cos \varphi+y_1,
         \sqrt{2I} \sin \varphi + y_2 \right).
   \end{align*}
   This system depends on one fast phase $ \varphi $ and does not depend
   analytically on $ I $ near $ I=0 $. This is not a minor problem since for
   many interesting problems coming from physics, a neighborhood of $ I=0 $
   may play the most important role. In the problem of semiclassical
   quantization, for examples it corresponds to the so-called low
   levels of Landau which are connected with Hall's conductivity. This problem
   is actually the main motivation of present paper.

    The procedure in \cite{Ne} is based upon subsequent change of
   variables, which corresponds to methods of KAM-theory.
   The aim of this paper is to show that one can choose a transformation
   such that the averaging procedure of \cite{Ne} is applicable in a
   neighborhood of $ I=0 $, and such that the passage from
   $ I > \varkappa > 0 $ to $ I=0 $ is uniform.

\section{ Formulation of the problem and the main result. }

   Consider the Hamiltonian
   \begin{align}
      H &= \mathcal{H}_0 \left( I \right)+
         \varepsilon g_0(q,p,y_1,y_2), & I &= \frac{q^2+p^2}{2},
     \label{IntrInitHam}
   \end{align}
   where $ 0 < \varepsilon < \varepsilon_0 $ is a small parameter, and
   $ \mathcal{H}_0$ and $ g_0 $ are real analytic functions
   in a complex $ \delta $-neighborhood of the domain
   $ D := D_{2n} \{ y_1,y_2 \} \times D_{2}^{\varkappa} \{ q,p \} $,
   $ D_{2n} \subset \mathbb{R}^{2n} $,
   $ D_{2}^{\varkappa} = \{ (q,p) \in \mathbb{R}^2 \, |
   \, I < \varkappa \} $. In $ D $ we assume the following conditions
   \begin{align*}
      | \mathcal{H}_0 | \le C, \qquad | g_0 | \le C, \qquad
         \left|\frac{\partial \mathcal{H}_0 }{\partial I}
     \right|  \ne 0.
   \end{align*}

   Following \cite{BM,Ne}, we show that for every integer $ m > 0 $ there
   exists a close to identity real analytic canonical transformation
   $ (q,p,y_1,y_2) \to (Q,P,z_1,z_2) $ defined by
   \begin{equation}
      \left\{
      \begin{aligned}
         q&= Q+\varepsilon Q^1(Q,P,z_1,z_2,\varepsilon), & \quad
        y_1 &= z_1 +\varepsilon Z_1^1(Q,P,z_1,z_2,\varepsilon), \\
         p&= P+\varepsilon P^1(Q,P,z_1,z_2,\varepsilon), & \quad
            y_2&= z_2 +\varepsilon Z_2^1(Q,P,z_1,z_2,\varepsilon), \\
      \end{aligned} \label{IntrCanTr} \right.
   \end{equation}
   where $ |Q^1| + | P^1 | + | Z_1^1 | + | Z_2^1 | \le C $, and such that
   the Hamiltonian (~\ref{IntrInitHam}~) transforms into
   \begin{align}
      H &= \mathcal{H}_m \left( \frac{Q^2+P^2}{2},z_1,z_2,
         \varepsilon \right) + \varepsilon g_m (Q,P,z_1,z_2,\varepsilon).
         \label{IntrFinalHam}
   \end{align}
   Precisely, we have the following theorem.

   {\bfseries Theorem 1. }
      Assume the conditions above in
      $ (Q,P,z_1,z_2) \in D+\frac{1}{2}\delta $.
      Then there exists some interval $ (0,\varepsilon_1 ] $, integer
      number $ r $ and real
      analytic canonical change of variables of the form(~\ref{IntrCanTr}~),
      which transforms Hamiltonian (~\ref{IntrInitHam}~) into
      (~\ref{IntrFinalHam}~) with exponentially small $ g $:
      \begin{align}
         | & g_{r} |  + | \nabla g_{r} | <
        c_2 \exp{(-\frac{1}{c_1 \varepsilon})}, &
            | & Q^1| + |Z_2^1|+|P^1|+|Z_1^1| < c_3, &
            | & \mathcal{H}_{r}- \mathcal{H}_0| < c_4 \varepsilon.
         \label{ThCondEst}
      \end{align}
      Here $ \varepsilon \in (0,\varepsilon_1] $ and
      $\varepsilon_1,r, c_i, i=1,2,3,4 $ depend on
      $ \varepsilon_0,\delta,C$ and $ \varkappa $.

  \section{ Auxiliary lemmas. }


   Let $ w \left(I,\mu \right)  $ and $ g (q,p,\mu) $ be an
   analytic function of $ I $ and $ (q,p) $ respectively and
   vector-parameter $ \mu $~, $ \frac{\partial w}{\partial I} \ne 0 $~,
   $ I= (q^2+p^2)/2 $. Denote
   $ \frac{\partial w}{\partial q} = q \frac{\partial w}{\partial I} $,
   $ \frac{\partial w}{\partial p} = p \frac{\partial w}{\partial I} $.
   Consider the equation~:
   \begin{align}
      & \frac{\partial w}{\partial p}\frac{\partial W}{\partial q} -
         \frac{ \partial w}{\partial q}\frac{\partial W}{\partial p} +
         g(q,p,\mu) = \bar g \left( \frac{q^2+p^2}{2},
         \mu \right), \label{LemmaInitEq} \\
      & \bar g (\frac{q^2+p^2}{2},\mu)= \left. 
        \int_{0}^{2\pi} g \left( q(\varphi,I),
        p(\varphi,I),\mu) \right. d\varphi 
	\right|_{\binom{\varphi=\varphi(q,p)}{I=I(q,p)}}, \\
      & g (q,p,\mu)=\bar g \left( \frac{q^2+p^2}{2},\mu \right) -
         \tilde g (q,p,\mu), \label{GammaBaTDef}
   \end{align}
   where $ \bar g ((q^2+p^2)/2,\mu) $ is a mean value of $ g (q,p,\mu) $
   with respect to $ \varphi $, and the rest part of $g $ is noted as 
   $ \tilde g (q,p,\mu) $ and it is convenient to take sign minus here.

  { \bfseries Lemma 1. }
  Equation (\ref{LemmaInitEq}) is solvable and has analytic solution
  on variables $ (q,p) $ and parameter $ \mu $. Function $ W(q,p,\mu) $,
  defined by formula~:
  \begin{align}
   W(q,p,\mu)  = \left. \frac{1}{\frac{\partial w}{\partial I}}
      \left( \frac{1}{2} \int\limits_0^{\varphi} \tilde g (q(\psi,I),
      p(\psi,I),\mu) d\psi + \frac{1}{2} \int\limits_{\pi}^{\varphi}
      \tilde g (q(\psi,I),p(\psi,I),\mu) d\psi \right)
      \right|_{\varphi=\varphi(q,p),I=I(q,p)}, \label{Lemma1CondSol}
  \end{align}
  is an analytic solution of equation (\ref{LemmaInitEq}).

  { \bfseries Remark. }
  General solution of equation (\ref{LemmaInitEq}) has a form~:
  \begin{align*}
      W(q,p,\mu) & = \left. \left( \frac{1}{\frac{\partial w}{\partial I}}
      \int\limits_0^{\varphi} \tilde g (q(\psi,I),
      p(\psi,I),\mu) d\psi + W_0(I,\mu) \right)
      \right|_{\varphi=\varphi(q,p),I=I(q,p)} ,
  \end{align*}
  where $ W_0(I,\mu) $ is a constant of integration.
   It follows for $ (q,p) $ variables from method \cite{Ne}.
   We show below, that in the neighborhood of point $ q=0, p=0 $
   definition of constant of integration, pointed in (~\ref{Lemma1CondSol}~),
   allow us to integrate equation (~\ref{LemmaInitEq}~),
   preserving analyticity of solutions.

   { \begin{center} \bfseries Proof. \end{center}}

  First make a canonical change of variables:
  \begin{align*}
     & q = \frac{u+i v}{\sqrt{2}}, &
        & p = \frac{v+i u}{\sqrt{2}}.
  \end{align*}

  Then equation (\ref{LemmaInitEq}) is :
  \begin{align}
     \left( \frac{\partial w}{\partial v}
        \frac{\partial W}{\partial u} -
        \frac{ \partial w}{\partial u}
    \frac{\partial W}{\partial v} \right)
         = \tilde g \left( \frac{u+i v}{\sqrt{2}},
    \frac{v+i u}{\sqrt{2}}, \mu \right). \label{LemmaUVEq}
  \end{align}

  Now we should expand $ \tilde g $ into Taylor series.
  The expansion has form:
  \begin{align*}
     \tilde g \left( \frac{u+i v}{\sqrt{2}},
    \frac{v+i u}{\sqrt{2}} , \mu \right) & =
    \sum_{k,l \in \mathbb{N}}
        \tilde g_{kl} (\mu) u^k v^l \frac{(k+l)!}{k! l!}.
  \end{align*}

  Using coordinates $ (\varphi,\rho) $: $ u=\rho e^{i \varphi} $,
  $ v= \rho e^{-i \varphi} $
  one can easily integrate equation (\ref{LemmaUVEq}):
  \begin{align}
     W &=\frac{\rho}{\frac{\partial w}{\partial \rho}}
        \left( \int\limits_0^{\varphi}
        \tilde g \left( \frac{u+i v}{\sqrt{2}},
    \frac{v+i u}{\sqrt{2}}, \mu \right) d\psi \right) +
     W_0 (\rho,\mu) = \notag \\
       &= \frac{\rho}{\frac{\partial w}{\partial \rho}} \left(
        \sum_{k,l \in \mathbb{N}} \tilde g_{kl}(\mu)
    \frac{u^k(\varphi,\rho) v^l(\varphi,\rho)}{i(k-l)}
    \frac{(k+l)!}{k!l!} - \sum_{k,l \in \mathbb{N}}
    \tilde g_{kl}(\mu) \frac{|uv|^{(k+l)}}{i(k-l)}
    \frac{(k+l)!}{k!l!} \right) +W_0(\rho,\mu),
    \label{LemmaWExpan}
  \end{align}
  where $ W_0(\rho,\mu) $ is a constant of integration.

  We can see that at bottom limit of integration we can have
  nonanalytic $ (~ \sim uv=\sqrt{2I}~)$ dependence on $ I $
  at the point $ (u=0,v=0) $. But we can choose
  $ W_0(\rho,\mu) $ in a such way that nonanalytic term
  disappears. Introduce $ W_0 $ as :
  \begin{align}
     W_0(\rho,\mu) &= \frac{\rho}{\frac{\partial w}{\partial \rho}}
        \left( \frac{1}{2} \int\limits_{\pi}^{\varphi}
    \tilde g (u(\psi,\rho),v(\psi,\rho),\mu) d\psi -
    \frac{1}{2} \int\limits_{0}^{\varphi}
    \tilde g (u(\psi,\rho),v(\psi,\rho),\mu) d\psi \right) = \notag \\
        & = \frac{\rho}{\frac{\partial w}{\partial \rho}} \left(
        \sum_{k,l \in \mathbb{N}} \tilde g_{kl}(\mu)
    \frac{|uv|^{(k+l)}}{i(k-l)} \frac{(k+l)!}{k!l!}
    \right). \label{LemmaW0Expan}
  \end{align}
  Hence (\ref{LemmaW0Expan}) 'kills' nonanalytic term
  in (\ref{LemmaWExpan}), so $ W $ is analytic function of $ I $ in the
  point $ I=0 $. It is easy to see that
  $ \frac{\partial w}{\rho \partial \rho}=\frac{\partial w}{\partial I}$.
  Function $ W(q,p,\mu) $ is analytic ones of parameter $\mu $ due to
  procedure of construction: we had no nonanalytic
  dependence on $ \mu $ on all steps of obtaining of solution.
  Hence, $ W $ is an analytic function of $ (q,p) $ and parameter
  $\mu$. Lemma is proved.


 \vspace{7mm}

  We use generating function to construct
  canonical transformation in the proof of main theorem.

  Let $ S(q,P,y_1,z_2,\varepsilon) $ be an analytic function
  of all variables $ (q,P,y_1,z_2) $ and small parameter
  $ \varepsilon \in [0,\varepsilon_0] $, $ (q,p,y_1,y_2) \in U $,
  $(Q,P,z_1,z_2) \in U-\delta $, $\delta >0$, $ U \subset
  \mathbb{R}^{2n+2} $. Domain $ U-\delta $ is the set of points from $ U $,
  which enter together with their $\delta $-neighborhood \cite{A2}. Consider
   a system of equations :
  \begin{equation}
     \left\{
     \begin{aligned}
        Q&=q+\varepsilon
           \frac{ \partial S(q,P,y_1,z_2,\varepsilon)}{\partial P}, & \quad
       z_1 &= y_1 + \varepsilon
       \frac{ \partial S(q,P,y_1,z_2,\varepsilon)}{\partial z_2}, \\
        p&=P+\varepsilon
           \frac{ \partial S(q,P,y_1,z_2,\varepsilon)}{\partial q}, & \quad
           y_2 &= z_2 + \varepsilon
       \frac{ \partial S(q,P,y_1,z_2,\varepsilon)}{\partial y_1}, \\
    \end{aligned}
    \right.  \label{Lemma2InitSys}
  \end{equation}
  and assume that $ \max \left\{ |\frac{\partial S}{\partial q}|,
  |\frac{\partial S}{\partial P}|, |\frac{\partial S}{\partial y_1}|,
  |\frac{\partial S}{\partial z_2}| \right\} < C_U $ in domain
  $ U  $.

  Following lemma establishes that introduced transformation is a
  change of variables and allows us to estimate corrections.

  { \bfseries Lemma 2. }
  If $ \varepsilon < \delta / (2 C_U (n+1) ) $, where $ n $ is a dimension
  of $ y_1 $, then the system
  (\ref{Lemma2InitSys}) is solvable and solution has
   the following form:
   \begin{equation}
      \left\{
      \begin{aligned}
         q&= Q+\varepsilon Q^1(Q,P,z_1,z_2,\varepsilon), & \quad
        y_1&= z_1 +\varepsilon Z_1^1(Q,P,z_1,z_2,\varepsilon), \\
         p&= P+\varepsilon P^1(Q,P,z_1,z_2,\varepsilon), & \quad
            y_2&= z_2 +\varepsilon Z_2^1(Q,P,z_1,z_2,\varepsilon), \\
      \end{aligned}
      \right.  \label{Lemma2Sol}
   \end{equation}
   where $ Q^1, P^1, Z_1^1, Z_2^1 $ are analytic functions in
   $ U-\delta \times [0,\varepsilon_0] $
   and $ \max \left\{ |Q^1|, |P^1|, |Z_1^1|,|Z_2^1| \right\} < C_U $
   in domain $ U-\delta $.

  { \begin{center} \bfseries Proof. \end{center}}

  Let us present system (\ref{Lemma2InitSys}) as~:
  \begin{equation}
     \left\{
     \begin{aligned}
        f_1&(q,p,y_1,y_2;Q,P,z_1,z_2,\varepsilon)=0, & \quad
       f_3&(q,p,y_1,y_2;Q,P,z_1,z_2,\varepsilon)=0, \\
        f_2&(q,p,y_1,y_2;Q,P,z_1,z_2,\varepsilon)=0, & \quad
       f_4&(q,p,y_1,y_2;Q,P,z_1,z_2,\varepsilon)=0, \\
     \end{aligned}
     \right. \label{Lemma2ImplFunSys}
  \end{equation}
  and $ F=(f_1,f_2,f_3,f_4) $.
  Due to the theorem of implicit function, the system
  (\ref{Lemma2ImplFunSys}) is solvable with respect to $ (Q,P,z_1,z_2) $,
  if $ \det|\frac{DF}{D(q,p,y_1,y_2)}| \ne 0 $, and its solution
  is analytic in $ U-\delta \times [0,\varepsilon_0] $.
  In our case
  \begin{align*}
     & \det \left| \frac{DF}{D(q,p,y_1,y_2)} \right| = \det \left|
     \begin{aligned}
        1 & +\varepsilon \frac{\partial^2 S}{\partial P \partial q} &
          &\varepsilon \frac{\partial^2 S}{\partial P \partial y_1} \\
     &\varepsilon \frac{\partial^2 S}{\partial z_2 \partial q} &
        1 + \varepsilon & \frac{\partial^2 S}{\partial z_2 \partial y_1} \\
     \end{aligned}
     \right|.
  \end{align*}
  $ S(q,P,y_1,z_2,\varepsilon) $ is an analytic function of
  $\varepsilon $. Using Cauchy estimations for analytic functions
  \cite{A}, we obtain
  \begin{align}
     \det \left| \frac{DF}{D(q,p,y_1,y_2)} \right|-1 <
     \sum_{i=1}^{n+1} \left( \varepsilon \frac{ C_U }{\delta} \right)^i
     \frac{(n+1)!}{(n+1-i)!} < \sum_{i=1}^{\infty} \left( \varepsilon
     \frac{ C_U }{\delta} (n+1) \right)^i < 1 \label{Lemma2DetIneq}.
  \end{align}
  Hence, if $ \varepsilon < \delta / \left( 2 C_U (n+1) \right) $  then
  (~\ref{Lemma2DetIneq}~) is valid on $ U $ always and our system is solvable.
  Uniqueness of the transformation is established similarly to \cite{A2}.
  The solutions should satisfy the equations
  (\ref{Lemma2InitSys}) with $ \varepsilon =0 $ also, so they can
  be represented as (\ref{Lemma2Sol}). Let us substitute the
  solutions in equations (\ref{Lemma2InitSys}).
  Then from equations for $ q $ and $ y_1 $ we have the following
  form for $ Q^1 $ and $ Z_1^1 $:
  \begin{align*}
     Q^1 & = -\frac{\partial S(Q+\varepsilon Q^1,P,z_1+\varepsilon Z_1^1,z_2,
        \varepsilon)}{\partial P}, &
     Z_1^1 & = -\frac{\partial S(Q+\varepsilon Q^1,P,z_1+\varepsilon Z_1^1,z_2,
        \varepsilon)}{\partial z_2}.
  \end{align*}
  Hence $ | Q^1 | < C_U $ and $ | Z_1^1 | < C_U $. For $ p $ and $ y_2 $
  variables we have the following formulae:
  \begin{align*}
     p = P + \varepsilon \frac{\partial S(Q+\varepsilon Q^1,P,z_1+\varepsilon Z_1^1,
        z_2,\varepsilon)}{\partial q}
       = P + \varepsilon \frac{\partial S(Q+\varepsilon Q^1,P,z_1+\varepsilon Z_1^1,
        z_2,\varepsilon)}{\partial Q}, \\
     y_2 = z_2 + \varepsilon \frac{\partial S(Q+\varepsilon Q^1,P,
        z_1+\varepsilon Z_1^1,z_2,\varepsilon)}{\partial y_1}
       = z_2 + \varepsilon \frac{\partial S(Q+\varepsilon Q^1,P,
        z_1+\varepsilon Z_1^1,z_2,\varepsilon)}{\partial z_1}.
   \end{align*}
   Therefore:
   \begin{align*}
      P^1 &=\frac{\partial S(Q+\varepsilon Q^1,P,z_1+\varepsilon Z_1^1,z_2,
         \varepsilon)}{\partial Q}, &
      Z_2^1 &=\frac{\partial S(Q+\varepsilon Q^1,P,z_1+\varepsilon Z_1^1,z_2,
         \varepsilon)}{\partial z_1},
    \end{align*}
   and again $ | P^1 | < C_U $, $  | Z_2^1 | < C_U $. Lemma is proved.


  \section{ Proof of the theorem. }

  We construct our change of variables as a composition of
  large number of consequently defined canonical transformations,
  giving dependence on $ (q^2+p^2)/2 $  of Hamiltonian more and
  more higher degree of $\varepsilon $.


  \subsection{ Procedure of consequently defined changes of
  variables. }

  Assume that Hamiltonian is obtained after $ i$ changes of
  variables has the following form~:
  \begin{align}
     H & =\mathcal{H}_i \left( \frac{q^2+p^2}{2},y_1,y_2,
        \varepsilon \right) + \varepsilon g_i(q,p,y_1,y_2,\varepsilon)
        \label{PTInitHam}, \\
     (q, & y_2,p,y_1) \in D_i, &  D_i &= D_1-2(i-1)K\varepsilon,
        \notag \\
     D&+\frac{\delta}{2}  \subset D_i \subset D+\delta, &
        D_1 &= D+\frac{3}{4}\delta. \notag
  \end{align}

  At $ i+1 $ step one has to find canonical infinitesimal (``almost identical'')
  change of variable (see \cite{Ne,A,A2}) $ (q,p,y_1,y_2) \to (Q,P,z_1,z_2) $
  \begin{equation}
     \left\{
     \begin{aligned}
     q&= Q+\varepsilon Q^1(Q,P,z_1,z_2,\varepsilon), & \quad
        y_1&= z_1 +\varepsilon Z_1^1(Q,P,z_1,z_2,\varepsilon), \\
     p&= P+\varepsilon P^1(Q,P,z_1,z_2,\varepsilon), & \quad
        y_2&= z_2 +\varepsilon Z_2^1(Q,P,z_1,z_2,\varepsilon), \\
     \end{aligned}
     \right. \label{PTCanSol}
  \end{equation}
  such that Hamiltonian takes a form~:
  \begin{align}
     & H =\mathcal{H}_{i+1} \left( \frac{Q^2+P^2}{2},z_1,z_2,
        \varepsilon \right)+
        \varepsilon g_{i+1} (Q,P,z_1,z_2,\varepsilon), \label{PTNextDef} \\
     & \mathcal{H}_{i+1} = \mathcal{H}_i + \varepsilon \bar g_i,
        \qquad \qquad g_{i+1}=O(\varepsilon^{i+1}), \notag
  \end{align}
  where $ \mathcal{H}_{i+1} $ contains terms of order $ i+1 $ of
  $ \varepsilon $, and the operation ``bar'' ´over $g_i$ is defined in sec. 3.

  It is possible and convenient to do by means of
  generating function (then one has the canonical property automatically)~:
  \begin{align*}
     S=S(q,P,y_1,z_2,\varepsilon)=qP+y_1z_2+
        \varepsilon S^1(q,P,y_1,z_2,\varepsilon).
  \end{align*}
  All other variables are defined in terms of $ (q,P,y_1,z_2)$ and
  $ S(q,P,y_1,z_2) $:
  \begin{equation}
     \left\{
     \begin{aligned}
        Q &=q+\varepsilon \frac{\partial S^1(q,P,y_1,z_2,\varepsilon)}
       {\partial P}, & \quad
       z_1 &=y_1+\varepsilon \frac{\partial S^1(q,P,y_1,z_2,\varepsilon)}
       {\partial z_2}, \\
        p &=P+\varepsilon \frac{\partial S^1(q,P,y_1,z_2,\varepsilon)}
       {\partial q}, & \quad
           y_2 &=z_2+\varepsilon \frac{\partial S^1(q,P,y_1,z_2,\varepsilon)}
       {\partial y_1}. \\
     \end{aligned}
     \right. \label{PTGFSol}
  \end{equation}

  Substitute (\ref{PTGFSol}) to (\ref{PTInitHam}) and (\ref{PTNextDef}),
  and equate Hamiltonians in ``mixed'' new-old´coordinates~:
  \begin{align*}
     & \mathcal{H}_i \left( \left( \frac{q^2+\left(P+
        \varepsilon \frac{\partial S^1}{\partial q}\right)^2}{2} \right),
    y_1, z_2+\varepsilon \frac{\partial S^1}{\partial y_1},
    \varepsilon \right)+
    \varepsilon g_i \left(q,P+\varepsilon \frac{\partial S^1}{\partial q},
        y_1,z_2+\varepsilon \frac{\partial S^1}{\partial y_1},\varepsilon
    \right) = \notag \\
     & \mathcal{H}_i \left( \left( \frac{ \left( q+
        \varepsilon \frac{\partial S^1}{\partial P} \right)^2 + P^2}{2}
    \right),
        y_1+\varepsilon \frac{\partial S^1}{\partial z_2},z_2,
    \varepsilon \right)+ \varepsilon \bar g_i \left( \frac{ \left( q+
        \varepsilon \frac{\partial S^1}{\partial P} \right)^2 + P^2}{2},
        y_1+\varepsilon \frac{\partial S^1}{\partial z_2},z_2,
    \varepsilon \right)+ \notag \\
      & +\varepsilon g_{i+1} \left(q+\varepsilon \frac{\partial S^1}{\partial P},P,
         y_1+\varepsilon \frac{\partial S^1}{\partial z_2},z_2,
     \varepsilon \right). \notag
   \end{align*}

   Expand them in Taylor series, write out the terms of the same degree 
   of $ \varepsilon $ and take into account, that 
   $ \partial \mathcal{H}_i/\partial y_1 = \partial \mathcal{H}_i/
   \partial z_2 = O (\varepsilon) $~:
   \begin{align}
      & \mathcal{H}_i \left( \frac{q^2+P^2}{2},y_1,z_2,\varepsilon \right)=
     \mathcal{H}_i \left( \frac{q^2+P^2}{2},y_1,z_2,\varepsilon \right), \\
      & \frac{\partial \mathcal{H}_i}{\partial P}
         \frac{\partial S^1}{\partial q} -
     \frac{\partial \mathcal{H}_i}{\partial q}
         \frac{\partial S^1}{\partial P} +g_i(q,P,y_1,z_2,\varepsilon) =
     \bar g_i \left( \frac{q^2+P^2}{2},y_1,z_2,\varepsilon \right).
     \label{PTGFEq}
   \end{align}
   Equation (\ref{PTGFEq}) is integrated with Lemma 1:
   $ S(q,P,y_1,z_2,\varepsilon)= W(q,P,\mu) $, where
   $ \mu=(y_1,z_2,\varepsilon) $. Function $ S $ is analytic function
   of all its variables and small parameter $\varepsilon $, so we can
   solve system (\ref{PTGFSol}), using Lemma 2. The solution is
   continuous and has continuous derivative, so it defines change of
   variables. We find necessary canonical transformation (\ref{PTCanSol})
   and after substitution into (\ref{PTInitHam}) we can obtain explicit
   expression for $ \mathcal{H}_{i+1} $ and $ \varepsilon g_{i+1} $.


   \subsection{ Estimations. }

   Let $ r $ steps have been done. Domain $ D_i $, where our
   Hamiltonian is considered after $i$ steps, defined as
   $ D_{i+1} = D_1-2(i-1)K\varepsilon $, where
   $ D_1= D+\frac{3}{4}\delta $, K -positive constant,
   defined below.

   At first step we can easily show, that if
   $ (q,p,y_1,y_2) \in D_1 $ and $ \varepsilon $ is quite small,
   then (\ref{PTCanSol}) is defined and the following conditions
   are realized:
   \begin{align*}
      \left| g_1 \right| + \left| \nabla g_1 \right| &< k_1\varepsilon, &
         \left| Q^1+P^1+Z_1^1+Z_2^1 \right| &< k_2, &
     \left| \mathcal{H}_{1}- \mathcal{H}_0 \right| &< k_3
     \varepsilon.
   \end{align*}
   $ \left\{ k_i \right\} $ are positive constants.
   Really, from (~\ref{PTNextDef}~) we immediately obtain that
   $ g_1= O ( \varepsilon^2 ) $. Using estimation procedure defined
   below all other inequalities are easily derived.

   Let us  accept inductive hypothesis that with $ i: 1\le i \le r $
   following estimations hold:
   \begin{align}
      \left| \mathcal{H}_i \right| &< 2 C, &  c &< \left|
         \frac{\partial \mathcal{H}_i}{\partial I} \right| < 2C, &
         \left| \nabla g_i \right| + \left| g_i \right| &< M_i, & \quad
         M_i&=\frac{k_1 \varepsilon}{2^{i-1}},
       \label{PTIndHyp}
   \end{align}
   where $ I=(q^2+p^2)/2 $.
    Now one should find $ \varepsilon_1 $ and $ K$ such that if
    $ 0<\varepsilon < \varepsilon_1 $,
    $ (Q,P,z_1,z_2) \in D_{r+1}=D_r-2K\varepsilon$, the system
    (\ref{PTCanSol}) is defined with $i=r$ and (\ref{PTIndHyp}) is
    satisfied with $ i=r+1 $.

    For  $ \frac{\partial S^1}{\partial \varphi} $ we have the following
    estimations from definition of $ \tilde g_i $ (~\ref{GammaBaTDef}~),
    equation (~\ref{PTGFEq}~) and form of solution (~\ref{Lemma1CondSol}~)~:
    \begin{align*}
       \left| \frac{\partial S^1}{\partial \varphi} \right| <
          \left| \tilde g_i \right| \le \left| g_i \right| < M_i.
    \end{align*}
    The estimations for $ \frac{\partial S^1}{\partial J} $,
    $ J=(q^2+P^2)/2 $ is obtained by the following way~:
    \begin{align}
       \frac{\partial S^1}{\partial J} = \frac{\partial}{\partial J}
        \left(  \frac{1}{\frac{\partial \mathcal{H}_i}{\partial J}} \left(
          \frac{1}{2}\int\limits_0^{\varphi} \tilde g (\psi,J) d\psi +
          \frac{1}{2}\int\limits_{\pi}^{\varphi} \tilde g (\psi,J) d\psi
           \right) \right) =\notag \\
       = -\frac{\partial^2 \mathcal{H}_i}{\partial J^2}
          \frac{S^1}{\frac{\partial \mathcal{H}_i}{\partial J}}+
          \frac{1}{\frac{\partial \mathcal{H}_i}{\partial J}}
      \frac{\partial}{\partial J} \left(
          \frac{1}{2} \int\limits_0^{\varphi} \tilde g (\psi,J) d\psi +
          \frac{1}{2} \int\limits_{\pi}^{\varphi} \tilde g (\psi,J) d\psi
      \right).  \label{PTdSdISum}
   \end{align}

   We can estimate $ \frac{\partial^2 \mathcal{H}_i}{\partial J^2} $ by~:
   \begin{align*}
      \left| \frac{\partial^2 \mathcal{H}_i}{\partial J^2} \right| <
         \left| \frac{\partial^2 \mathcal{H}_0}{\partial J^2} \right| +
     \varepsilon \left| \frac{\partial^2 \bar g_0}{\partial J^2} \right| +
     \varepsilon  \sum_{j=1}^i \left|
        \frac{\partial^2 \bar g_j}{\partial J^2} \right|.
   \end{align*}
   Using initial data we can state that
   $ | \frac{\partial^2 \mathcal{H}_0}{\partial J^2} | < m_1 $ and
   $ | \frac{\partial^2 \bar g_0}{\partial J^2} | < k_4 $, where
   $ \{m_i\} $ are positive constant and doesn't depend on step $ i $.
   For $ | \frac{\partial^2 \bar g_j}{\partial J^2} |  $
   we use Cauchy estimations for analytic functions \cite{A}~:
   \begin{align*}
      \left| \frac{\partial^2 \bar g_i}{\partial J^2} \right| &<
         \frac{M_i}{K \varepsilon}, & \sum_{j=1}^i \left|
     \frac{\partial^2 \bar g_i}{\partial J^2} \right| &<
     \sum_{j=1}^i \frac{k_1 \varepsilon}{2^{j-1} K \varepsilon } <
     \frac{2k_1}{K}.
   \end{align*}
   Summing up, we obtain:
   \begin{align}
      \left| \frac{\partial^2 \mathcal{H}_i}{\partial J^2} \right| < m_1+
      k_4\varepsilon+ \frac{2k_1}{K}\varepsilon < k_5. \label{PTd2GdI2Est}
   \end{align}
   For $ \frac{\partial g_i}{\partial J} $ we have the following
   estimations from (\ref{PTIndHyp}):
   \begin{align*}
      \left| \frac{\partial g_i}{\partial J} \right| =
         \left| \frac{\partial \rho}{\partial J} \right|
            \left| \frac{\partial g_i}{\partial \rho} \right| <
     \frac{1}{\rho} \left( \left| \frac{\partial g_i}{\partial P} \right|+
        \left| \frac{\partial g_i}{\partial q} \right| \right) <
     \frac{1}{\rho} 2 M_i,
   \end{align*}
   where $ \rho^2=J $.

   Then integral in (\ref{PTdSdISum}) can be estimated as:
   \begin{align*}
      \left| \frac{1}{\frac{\partial \mathcal{H}_i}{\partial J}} \left(
         \frac{1}{2} \int\limits_0^{\varphi}
     \frac{\partial \tilde g(\psi,J)}{\partial J} d\psi +
     \frac{1}{2} \int\limits_{\pi}^{\varphi}
     \frac{\partial \tilde g (\psi,J)}{\partial J} d\psi
     \right) \right| < m_2 (\frac{4 \pi M_i}{\rho}) < \frac{m_3 M_i}{\rho}.
   \end{align*}
   Then expression (\ref{PTdSdISum}) is estimated:
   \begin{align*}
      \left| \frac{\partial S^1}{\partial J} \right| <
         (k_5 +\frac{m_3}{\rho})M_i.
   \end{align*}

   Now we have the estimations for $ \frac{\partial S^1}{\partial \varphi} $
   and $ \frac{\partial S^1}{\partial J} $ and can obtain estimations for
   $ \frac{\partial S^1}{\partial q} $ and
   $ \frac{\partial S^1}{\partial P} $.
   \begin{align}
      \left| \frac{\partial S^1}{\partial P} \right| &= \left|
         \frac{\partial J}{\partial P}\frac{\partial S^1}{\partial J}+
     \frac{\partial \varphi}{\partial P}\frac{\partial S^1}{\partial \varphi}
     \right| < \rho (k_5 +\frac{m_3}{\rho})M_i + \left| \frac{q}{q^2+P^2} \right|
         \left| g_i \right| < \notag \\
      &< (\rho k_5 +m_3)M_i + \left| \frac{q}{q^2+P^2} \right|
            \left| \frac{\partial g_i}{\partial q} q+
         \frac{\partial g_i}{\partial P} P \right| < \notag \\
      &< (\rho k_5 +m_3)M_i + 2 M_i < m_4 M_i. \label{PTdSdPEst}
   \end{align}
   For $ \frac{\partial S^1}{\partial q} $ estimations obtain the same way.
   \begin{align}
      \left| \frac{\partial S^1}{\partial q} \right| < m_5 M_i. \label{PTdSdqEst}
   \end{align}
   Estimations for the
   $ \frac{\partial S^1}{\partial y_1} $ and $ \frac{\partial S^1}{\partial z_2} $
   are obtained easily by differentiation of $ \tilde g_i $ in the
   definition of solution for $ S $~:
   \begin{align}
      \left| \frac{\partial S^1}{\partial y_1} \right| &< m_6 M_i, &
         \left| \frac{\partial S^1}{\partial z_2} \right| &< m_7 M_i.
     \label{PTdSdyzEst}
   \end{align}
   Now we can define $ k_6= \max ( m_4,m_5,m_6,m_7) $.
   Function $ g_{i+1} $ can be defined by~:
   \begin{align*}
      \left| g_{i+1} \right| = & \left|
         \frac{\varepsilon}{2} \left(
     \frac{\partial^2 \mathcal{H}_i}{\partial P^2} \right)_{\theta}
     \left( \frac{\partial S^1}{\partial q} \right)^2 -
     \frac{\varepsilon}{2} \left(
     \frac{\partial^2 \mathcal{H}_i}{\partial q^2} \right)_{\theta}
     \left( \frac{\partial S^1}{\partial P} \right)^2  +
         \left( \frac{\partial \mathcal{H}_i}{\partial z_2} \right)_{\theta}
     \frac{\partial S^1}{\partial y_1} -
         \left( \frac{\partial \mathcal{H}_i}{\partial y_1} \right)_{\theta}
         \frac{\partial S^1}{\partial z_2} + \notag \right. \\
      & \left. +\varepsilon \left( \frac{\partial g_i}{\partial P}
         \right)_{\theta} \frac{\partial S^1}{\partial q}-
         \varepsilon \left( \frac{\partial \bar g_i}{\partial q} \right)_{\theta}   \frac{\partial S^1}{\partial P}+
         \varepsilon \left( \frac{\partial g_i}{\partial z_2} \right)_{\theta}
     \frac{\partial S^1}{\partial y_1} - \varepsilon \left(
     \frac{\partial \bar g_i}{\partial y_1} \right)_{\theta}
     \frac{\partial S^1}{\partial z_2} \right| ,
   \end{align*}
   where $ \left(\right)_{\theta} $ means the derivative in the middle point.
   Thus, $ g_{i+1} $ estimates using (\ref{PTIndHyp}),
   (\ref{PTd2GdI2Est})--(\ref{PTdSdyzEst})~:
   \begin{align*}
      \left| g_{i+1} \right| <
         \left| 2\frac{\varepsilon k_5 k_6 M_i^2}{2} +
     2 \varepsilon \left( k_4+\frac{2k_1}{K} \right) M_i +
     4 \varepsilon k_6 M_i^2 \right| < k_7 M_i. \varepsilon
   \end{align*}
   Using Cauchy estimations \cite{A} we obtain estimations for
   $ \nabla g_{i+1} $:
   \begin{align*}
      \left| \nabla g_{i+1} \right| &< \left|
          4\frac{\varepsilon k_5 k_6 M_i^2}{2 K \varepsilon} + 2 \varepsilon
         \left( k_4+\frac{2k_1}{K} \right) \frac{M_i}{K \varepsilon} +
      8 \varepsilon k_6 \frac{M_i^2}{K \varepsilon} \right| <
          k_8 M_i \left( \varepsilon+\frac{1}{K} \right).
   \end{align*}
   Choosing $ K $ quite large and $ \varepsilon $ quite small we
   obtain $ k_7 \varepsilon < \frac{1}{4} $,
    $ k_8 \left( \varepsilon + \frac{1}{K} \right) < \frac{1}{4} $.
   Then
   \begin{align*}
      \left| g_{i+1} \right| + \left| \nabla g_{i+1} \right| <
         \frac{M_i}{2}=M_{i+1},
   \end{align*}
   and other inductive inequalities (\ref{PTIndHyp}) are valid for $ i=r+1 $.
   Therefore, we can do necessary changes of
   variables with chosen $ K,\varepsilon $ till $ D_r $ is not empty. After
   $ r=(\frac{\delta}{4 K \varepsilon} ) > k_9/\varepsilon $ changes we have
   obtained:
   \begin{align*}
       \left| g_r \right|+\left| \nabla g_r \right| <
       \frac{k_1 \varepsilon}{2^{r-1}} < c_2 exp(-\frac{1}{c_1 \varepsilon}),
   \end{align*}
   and other inequalities (\ref{ThCondEst}) will be satisfied.

   Theorem is proved.

\section{Acknowledgements.}

   Authors express their heartfelt gratitude to V.Geyler,
   K.Pankrashkin, A.Shafarevich and to A.Neishtadt especially for
   useful discussions.

   The work is partially supported by the project
   DFG-RAS 436 RUS113/572 and FSB 282.

\end{document}